\documentclass[12pt,leqno]{article}
\usepackage{amsfonts}
\usepackage{graphicx}
\newcommand{\EE}{{\rm I\kern-2pt E}}
\newcommand{\RR}{{\rm I\kern-2pt R}}
\newcommand{\DD}{{\rm I\kern-2pt D}}
\newcommand{\PP}{{\rm I\kern-2pt P}}
\newcommand{\NN}{{\rm I\kern-2pt N}}
\newcommand{\dd}{{\rm \kern 3pt I\kern-9pt d}}

\title{SOME THOUGHTS UPON AXIOMATIZED LANGUAGES WITH EXTENSION TOOLS:\\
A Focus on Probability Theory and Error Calculus with Dirichlet Forms}
\author{Nicolas Bouleau}
\date{Lecture at the Catalan Studies Institute, 4th April 2003}
\begin{document}
\maketitle
\vspace{1cm}

This year marks the centenary of the birth of Kolmogorov. It is a pleasure for me to acknowledge
 this occasion by giving a lecture in connection with his life's work. My purpose herein is  certainly not to present
 a whole historical study of Kolmogorov's output, but rather provide some remarks on specific mathematical topics
 in which he played an active role. As you know, Kolmogorov produced some eight hundred publications
 encompassing all the main fields of mathematics: functional analysis, ergodic theory, turbulence, probability
theory and statistics, and logic. 
\begin{center}
\includegraphics[width=2in]{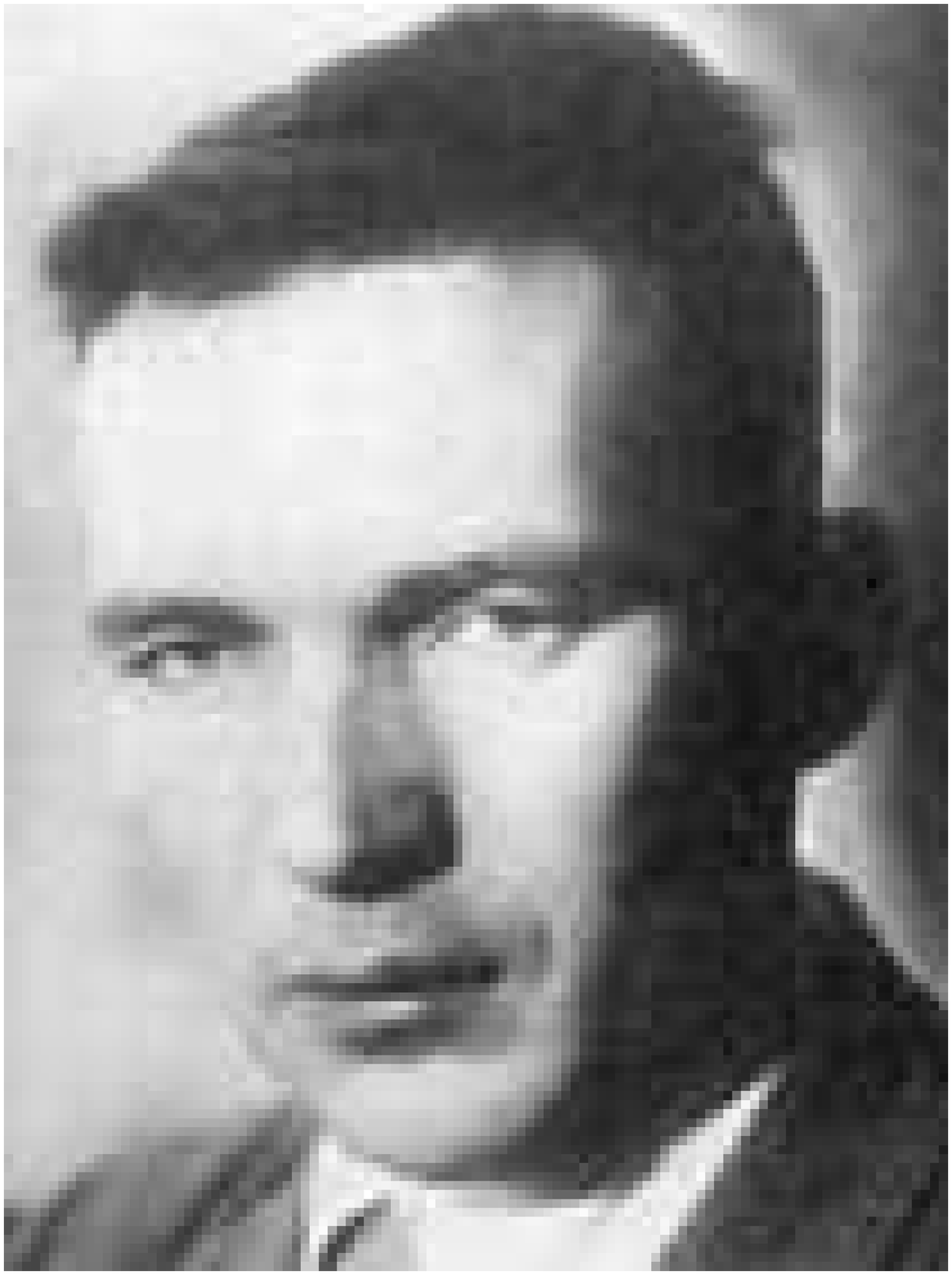}
\quad\quad\includegraphics[width=1.8in]{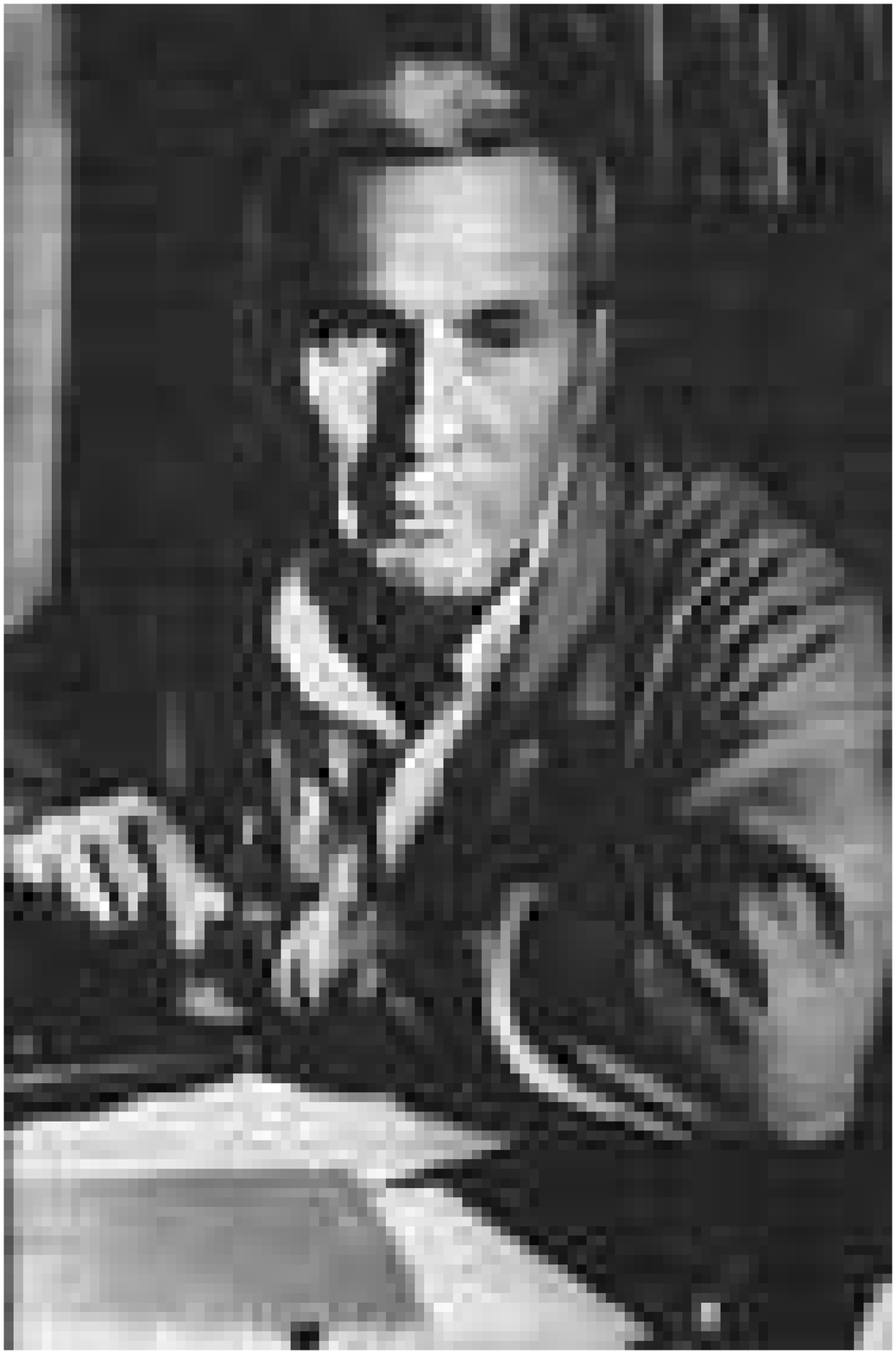}\\
Andre\" i Nikola\" ievitch Kolmogorov
\end{center}

He even delivred five seminal papers in the restricted domain of probability and stochastic
processes foundations between 1931 to 1936,  which make him one of the founders of the
 theory of continuous-time Markov processes or diffusions\footnote{These significant articles are the following:

- \"{U}ber die analytischen Methoden in der Wahrscheinlichkeitsrecgnung, 1931,

- Beitrage zur Masstheorie, 1933,

- Zur Theorie der stetigen zuf\"{a}lligen Prozessen, 1933,

- Grundbegriffe der Wahrscheinlichkeitsrechnung, 1993,

- Zur Theorie der Markoffschen Ketten, 1936.
}.
The subject I would like to discuss pertains to his famous "Grundbegriffe der
 Wahrscheinlichkeitsrechnung", which partially lies beyond his main body of mathematical work, 
 in some respects it serves as  {\it a manifesto} for how to tackle probability and probabilistic problems 
within the field of mathematics. I will be providing some remarks on axiomatized languages that display the cases
of both probability theory and of error calculus with Dirichlet forms. Based on these two examples, my aim
 is to emphasize
 the importance, in order for a language to be useful, of having an extension tool readily available.\\

\noindent I. {\large A brief history of random sequences theory}\\

In order to draw a comparison with Kolmogorov's axiomatic theory  , it is helpful to explain what 
 the "theory of random sequences" has become during the twentieth century. It did indeed serve an alternative
 way for incorporating probability into mathematics. Its purpose has been to describe a sequence of independent samples 
of a given
 quantity. In the simplest case, the theory pertains to samples of a random integer or even a single digit, so as to model the fair
game of heads and tails, in the one-half / one-half perfectly symmetric case\footnote{This section is inspired
 by the very interesting study conducted by Claude Dellacherie entitled ``Nombres au hasard de Borel \`a Martin L\"of"
{\it Gazette des Math\'ematiciens} n$^0 11$, 1978}.\\

\noindent I.1 {\it The normal numbers of Borel (1909)}\\

It is now easy, and Emile Borel was already able to make the proof  in 1909, that if we represent 
a real number over the unit interval $[0,1]$ by its binary expansion
$$(a_0, a_1, \ldots)\in\{0,1\}^{\NN}\quad\quad\longleftrightarrow\quad\quad x=\sum_{n=0}^\infty \frac{a_n}{2^{n+1}}\in[0,1]$$
to the independent one-half / one-half distribution of the digits corresponds the Lebesgue measure
on the interval $[0,1]$.
\begin{center}
\includegraphics[width=1.9in]{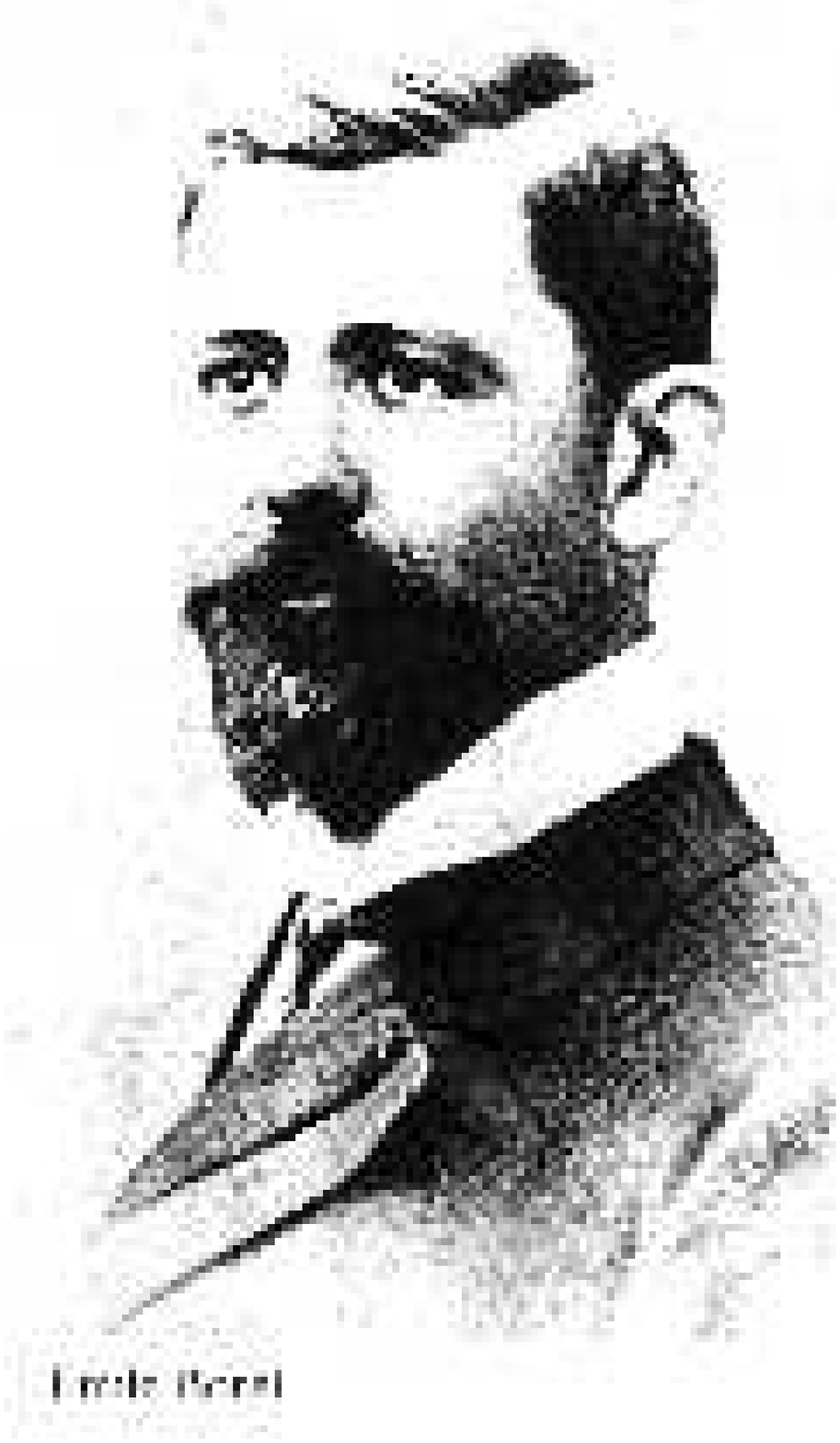}\\
Emile Borel\end{center}

As a consequence, for almost every real number $x\in[0,1]$, the asymptotic frequency
 of any finite sequence is $\frac{1}{2}$ to the power  of the sequence length. A real number fulfilling
 this property is said to be normal in the sense of Borel.

Now, proving that almost all real numbers are normal is just one step, another would be to exhibit such a number !
For the number $\pi$ determining whether it is normal or not  constitute a  famous  unsolved conjecture. Borel actually forwarded
 an effective, albeit sophisticated, construction of a normal number.

In 1933 however, Champernowne showed that the sequence obtained by writing the integer successively
 in dyadic representation is normal in the sense of Borel:
$$0\, 1\, 10\, 11\, 100\, 101\, 110\, 111\, 1000\, 1001\, 1010\,
 1011\, 1100\, 1101\, 1110\, 1111\, 10000\, ...$$
This clearly displays  that the concept of a normal number does not  capture the idea of random sequence very well.

Already back in 1919, Von Mises had proposed an improvement toward the definition a random sequence, 
by means of a new concept of ``collective"\footnote{``Grundlagen der Wahrscheinlichkeitsrechnung" {\it Math. Zeitung}
 5, 52-99, 1919.} which sought to describe a typical game of heads and tails. The idea is to ask for more
 than asymptotic averages and to think of a player gambling only at some random times depending on the 
evolution of the game : a sequence of digits is a ``collective" if it satisfies the law of large numbers and 
if any subsequence obtained by a non-anticipative selection rule satisfies also the law of large numbers. This
 interesting approach, which portends the notion  of ``stopping time", does nevertheless have the disadvantage of being difficult
  to apply in practical terms. A. Wald, one of the founders of statistics and decision theory, proposed in 1937 the more
 precise notion of ``collective relatively to a family of rules".
\begin{center}
\includegraphics[width=1.5in]{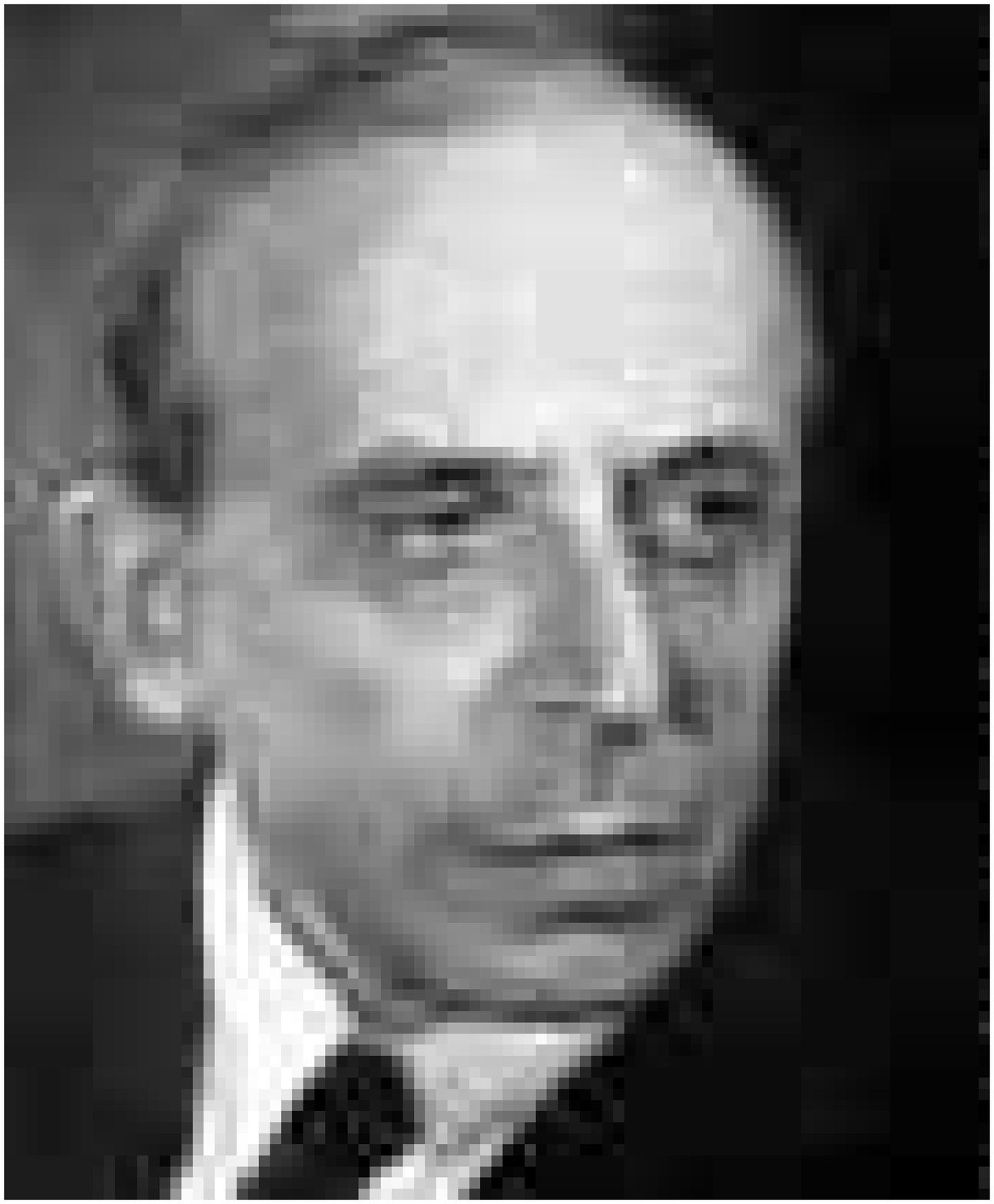}\quad\quad
\includegraphics[width=1.5in]{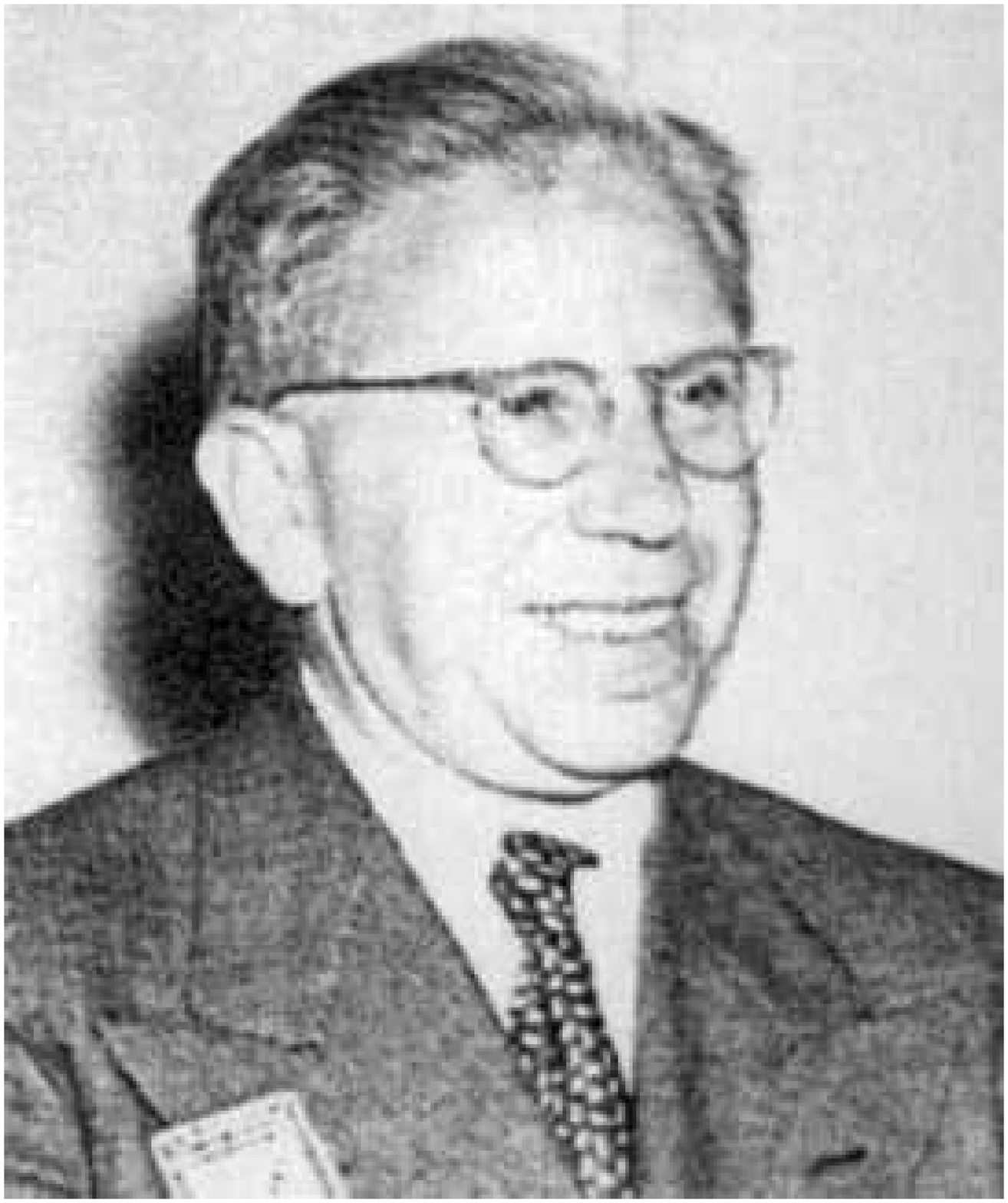}
\quad\quad
\includegraphics[width=1.5in]{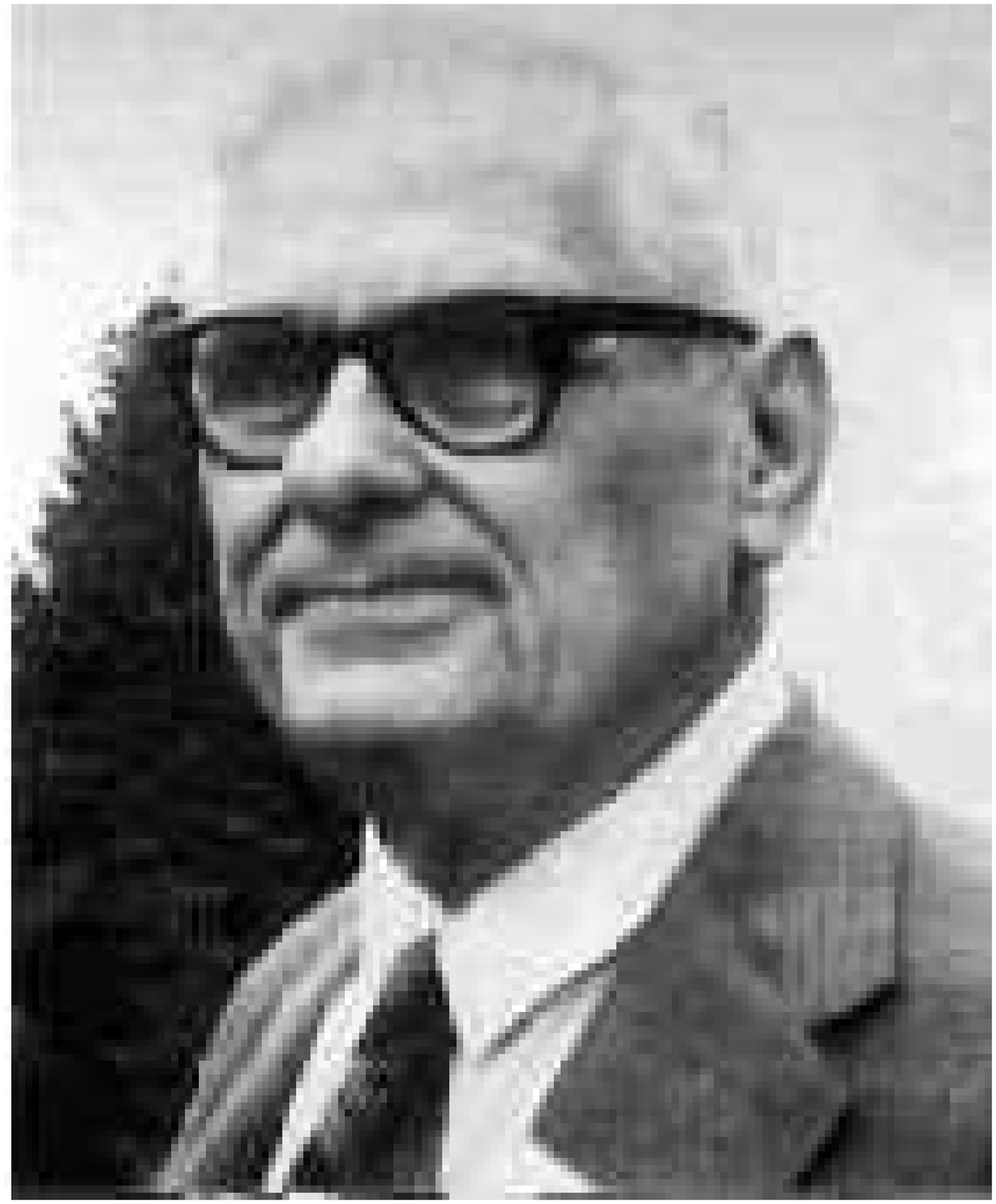}\\
R. Von Mises \qquad\qquad\qquad\quad A. Wald \qquad\qquad\qquad\quad A. Church\end{center}
Yet it would take the famous logician A. Church
 in 1940, with the first contribution from the field of logic into the debate, to propose an ``absolute notion of collective" 
that uses the set of all {\it effective} non-anticipative rules as regards recursive functions theory. It thus appeared
 that the goal has been achieved by applying this new theory of effectiveness stemming from the recent 
works of the logicians in the 1930's (G\"odel, Turing, Church).

Over this same period however, just prior to the Second World War, unsuspected new difficulties arose 
 concerning the notion of ``collective". In his work {\it Etude critique de la notion de collectif} (1939), Jean Ville
 demonstrated that random sequences  possess some probabilistic properties that a ``collective" may not always fulfill. 
A``collective"
does not  generally feature the right magnitude of fluctuations. In his argument Jean Ville uses the modern 
concept of mathematical {\it martingale} whose properties would be improved  by J. L. Doob in particular during the 1950's. By
 transfering the term {\it martingale} from gambling to mathematics Ville added a spark to this notion
 and likely contributed to its subsequent importance.

We would have to wait until the 1960's to obtain a satisfactory answer to the question of random sequence. This answer came
 from mathematical logic and is owed to Martin L\" of\footnote{``The definition of a random sequence" 
{\it Information and control} 9, 602-619, (1966).}. Roughly speaking, a random sequence 
 successfully passes all effective statistical randomness  tests. For a real number in $[0,1]$, being
 random in the sense of Martin L\" of signifies that it does not belong to any effective Lebesgue negligible
 set in $[0,1]$. Such a number cannot be given by an algorithm, it is random in the sense of Church yet 
avoids  Ville's critiques.

Although quite fascinating, the theory of random sequences remained useless for probabilists. The
 outstanding twentieth century development of probability theory, which began as a subsidiary
field and became one of the primary domains of applied and even pure mathematics, is based on another approach : 
the construction of a {\it language} for handling probabilistic calculations.

We would like to examine the reason behind this language's fruitfulness.\\

\noindent II. {\large Axiomatization of Kolmogorov and $\sigma$-additivity}\\

The paper entitled {\it Grundbegriffe der Wahrscheinlichkeitsrechnung} is an appeal to include probabilistic calculus
 into measure theory. Kolmogorov does not presume this idea is new, instead,  he cites several authors
 who have already applied Lebesgue measure theory for probabilistic investigations, in particular Borel, 
Fr\'echet, Steinhaus, L\'evy. He did proposes however new arguments, which proved to be highly valuable for
 subsequent research : the construction of probabilities on infinite dimensional spaces and the 
definition of conditional laws and conditional expectations using the Radon-Nikodym theorem.
\begin{center}
\includegraphics[width=1.7in]{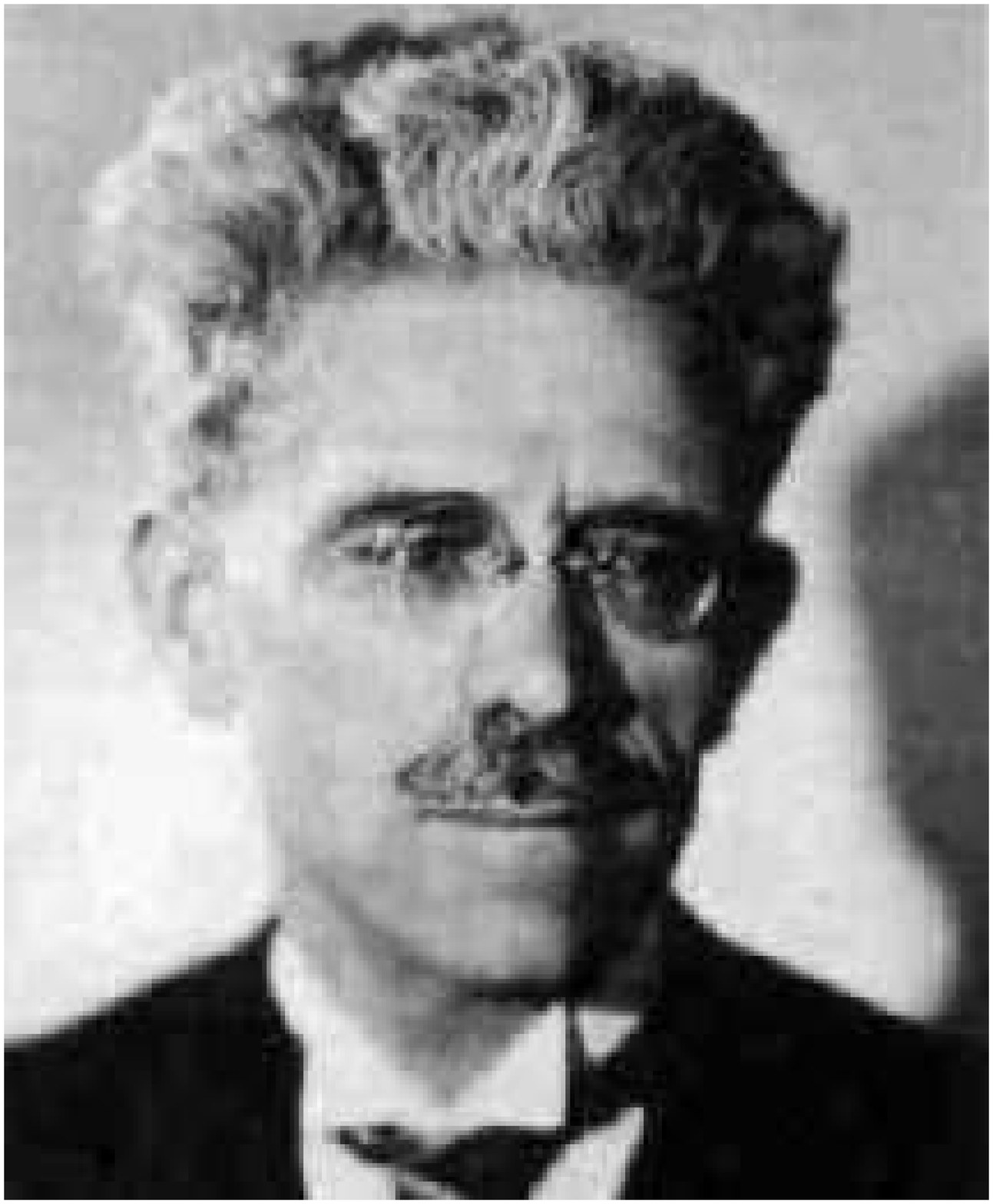}\quad\quad
\includegraphics[width=1.7in]{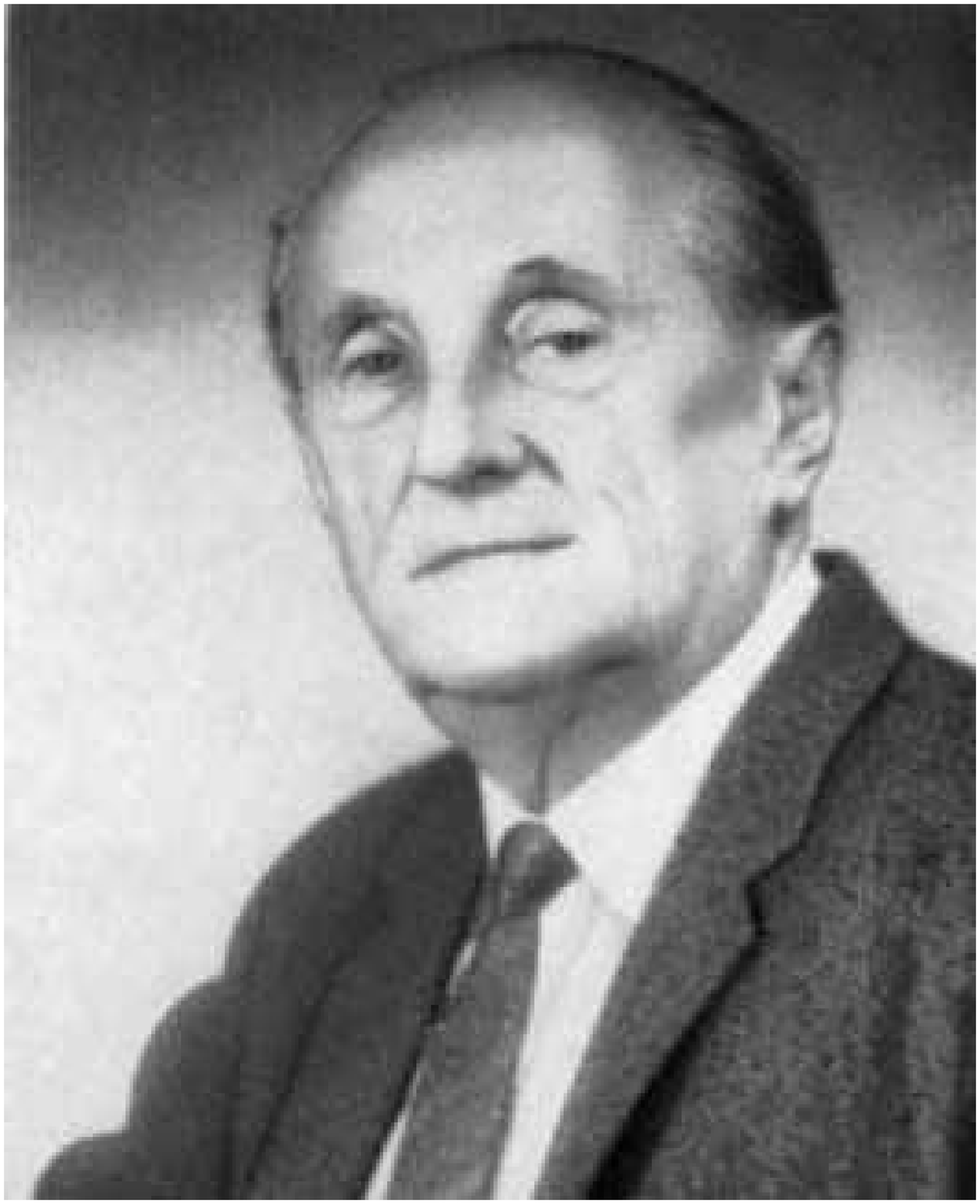}\quad\quad
\includegraphics[width=1.4in]{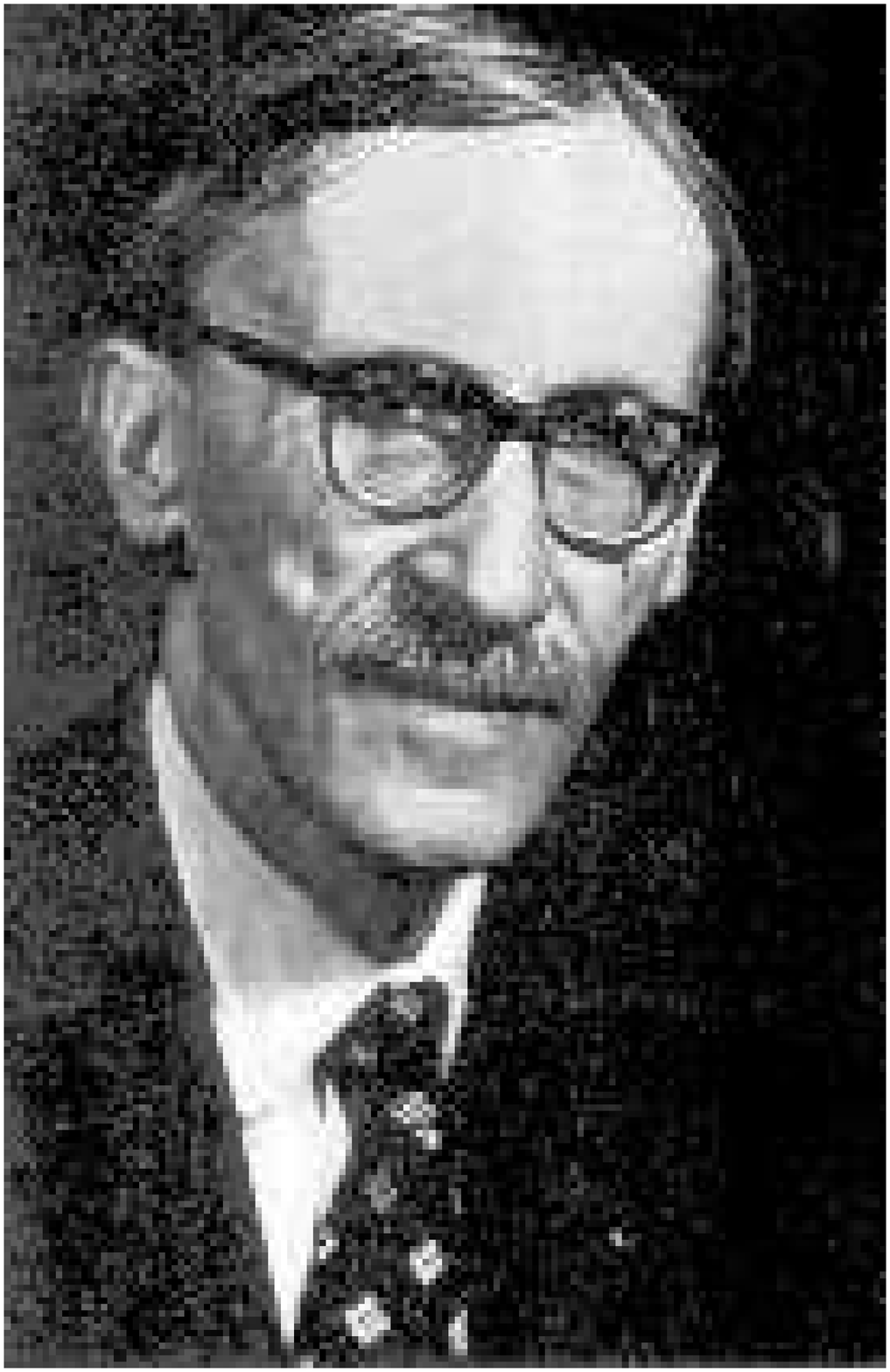}\\
M. Fr\'echet \qquad\qquad\quad\quad H. Steinhaus \qquad\qquad\quad\quad P. L\'evy\end{center}

He did not consider axiomatization as a pure formal system, but rather as a {\it language} that makes sense and that allows 
conducting thought and reasonning. In remarking that ``every axiomatic theory admits, as is well known,
an unlimited number of concrete interpretations"\footnote{We are in 1933 here and the works of L\" owenheim and 
Skolem (1915-1920) are already known, which prove  the existence of a countable model for any consistent theory.}, he emphasizes  
the intuitive interpretation of his axiomatization. He went on to display a dictionary between random events and sets:\\

\begin{tabular}{ll}
\multicolumn{1}{c}{\it Theory of sets}&\multicolumn{1}{c}{\it Random events}\\
1. $A$ and $B$ do not intersect, i.e. $AB=0$ & 1. Events $A$ and $B$ are incompatible\\
2. $AB\ldots N=0$ & 2. Events $A,B,\ldots, N$ are incompatible\\
3. $AB\ldots N=X$ & 3. Event $X$ is defined as the simultaneous \\
& occurrence of events $A,B,\ldots, N$\\
4. $A\cup B\cup\ldots\cup N=X$ & 4. Event $X$ is defined as the occurence \\
& of at least one of the events $A,B,\ldots, N$\\
5. The complementary set $A^c$ & 5. The non-occurence of event $A$\\
6. $A=0$ & 6. Event $A$ is impossible\\
7. $A=E$ & 7. Event $A$ must occur\\
8. Disjoint decomposition of $E$ & 8. Possible results $A_1, A_2, \ldots, A_n$\\
$A_1+A_2+\cdots+A_n=E$ & of an experiment\\
9. $B$ is a subset of $A$ & 9. From the occurence of event $B$\\
$B\subset A$ & follows the inevitable occurence of $A$
\end{tabular}\\

 For the axioms, the five first ones are elementary:
{\it Let ${\cal F}$ a set of subsets of a set $E$.

1. $\cal F$ is a field of sets

2. $\cal F$ contains the set $E$

3. To each set $A$ in $\cal F$ is assigned a non negative real number $P(A)$, called the probability of event $A$

4. $P(E)$ equals $1$

5. If $A$ and $B$ have no elements in common, then $P(A+B)=P(A)+P(A)$}\\

\noindent Kolmogorov underscores the importance of the sixth axiom : ``In all future investigations we shall assume that
 besides 
axioms 1 through 5,  another axiom holds true as well :

{\it 6. For a decreasing sequence of events 
$$A_1\supset A_2\supset\cdots\supset A_n\supset\cdots$$
in $\cal F$ for which $\cap_n A_n=0$ the following relation holds $\lim_n P(A_n)=0$".}\\

This axiom of $\sigma$-additivity implies the probability $P$ to be a measure in the sense of Lebesgue and Borel, which then
 embeds probability theory into measure theory :

\begin{center}
\begin{tabular}{ccc}
probability & $\longleftrightarrow$ & measure\\
event & $\longleftrightarrow$ & measurable set\\
random variable & $\longleftrightarrow$ & measurable function\\
expectation & $\longleftrightarrow$ & integral\\
independence & $\longleftrightarrow$ & product of measurable spaces\\
conditional expectation & $\longleftrightarrow$ & Radon-Nikodym derivative
\end{tabular}
\end{center}

Let's remark that as late as 1938, the philosopher Karl Popper, whose main education stemmed from the field of
psychology, was not convinced of the interest in placing probability theory within the framework of measure 
theory. Even in 1955, he still seemed proud to emphasize that a theory with only the first five axioms is more general. He
wrote ``Kolmogorov's system can be taken, however, as one of the interpretation of mine"\footnote{K. Popper, {\it The
logic of Scientific Discovery}, Hutchinson, 1972, p319.}.

\begin{center}
\includegraphics[width=1.5in]{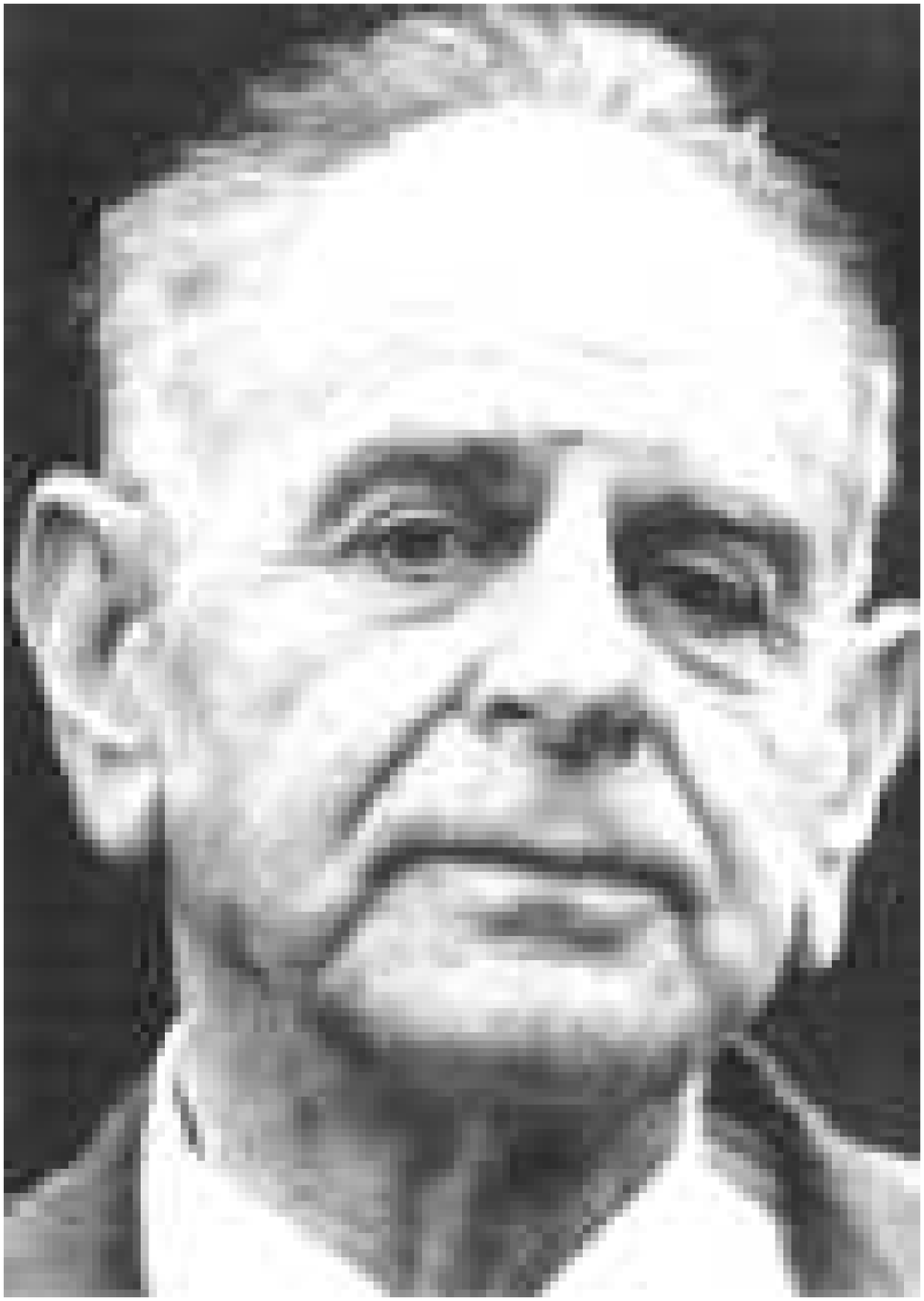}\\
Karl Popper
\end{center}

We know clearly now, thanks to the development of stochastic analysis over the twentieth
 century, that $\sigma$-additivity is the key tool making this language expansive. It allows defining the probability of 
events or the expectation of functions that are not given by simple closed formulae, but rather by limits. This fact is of absolutely
prime importance since several mathematical objects are defined by limits and {\it the methods for
 defining these converging sequences of  objects are not a priori restricted}.

This paves the way to the study of stochastic processes : if we know the probabilistic properties of a finite number of 
coordinates $X_n$ on a product space, without the $\sigma$-additivity we cannot conclude anything about functions depending
upon an infinite number of  $X_n$'s.

Thanks to  $\sigma$-additivity, connections with functional analysis may be developed, thereby giving rise to probabilistic
interpretations. For example, potential theory is connected with Markov processes theory and 
 martingales theory. Let's recall that J. L. Doob proved his extension of Fatou's lemma at the boundary from conical limits
to non-tangential limits, first using a probabilistic argument and then, one year later, by means of a purely 
analytical approach.\\

\noindent III. {\large Error calculus with Dirichlet forms}\\

I would  now like to present a more recent theory, in some repect a ``cousin" to probability theory, which also possesses  a means
 of 
extension providing it with remarkable power and fruitfulness. I have in mind the theory of Dirichlet forms with its interpretation
in terms of errors. I shall begin with the ideas of Gauss about errors which are the elementary bases
 of the theory.\\

\noindent III.1. {\it Gauss formulae for the propagation of errors}\\

The ideas of Gauss were forwarded at the beginning of the XIXth century, at a time when several mathematicians were 
concerned with  measurements errors, especially in the field of celestial mechanics. First of all, Legendre ({\it Nouvelles m\'ethodes
pour la dŽtermination des orbites des plan\`etes}, 1805) proposed the least squares principle   to choose the 
best value of a quantity obtained by several different measures. 
\begin{center}
\includegraphics[width=1.5in]{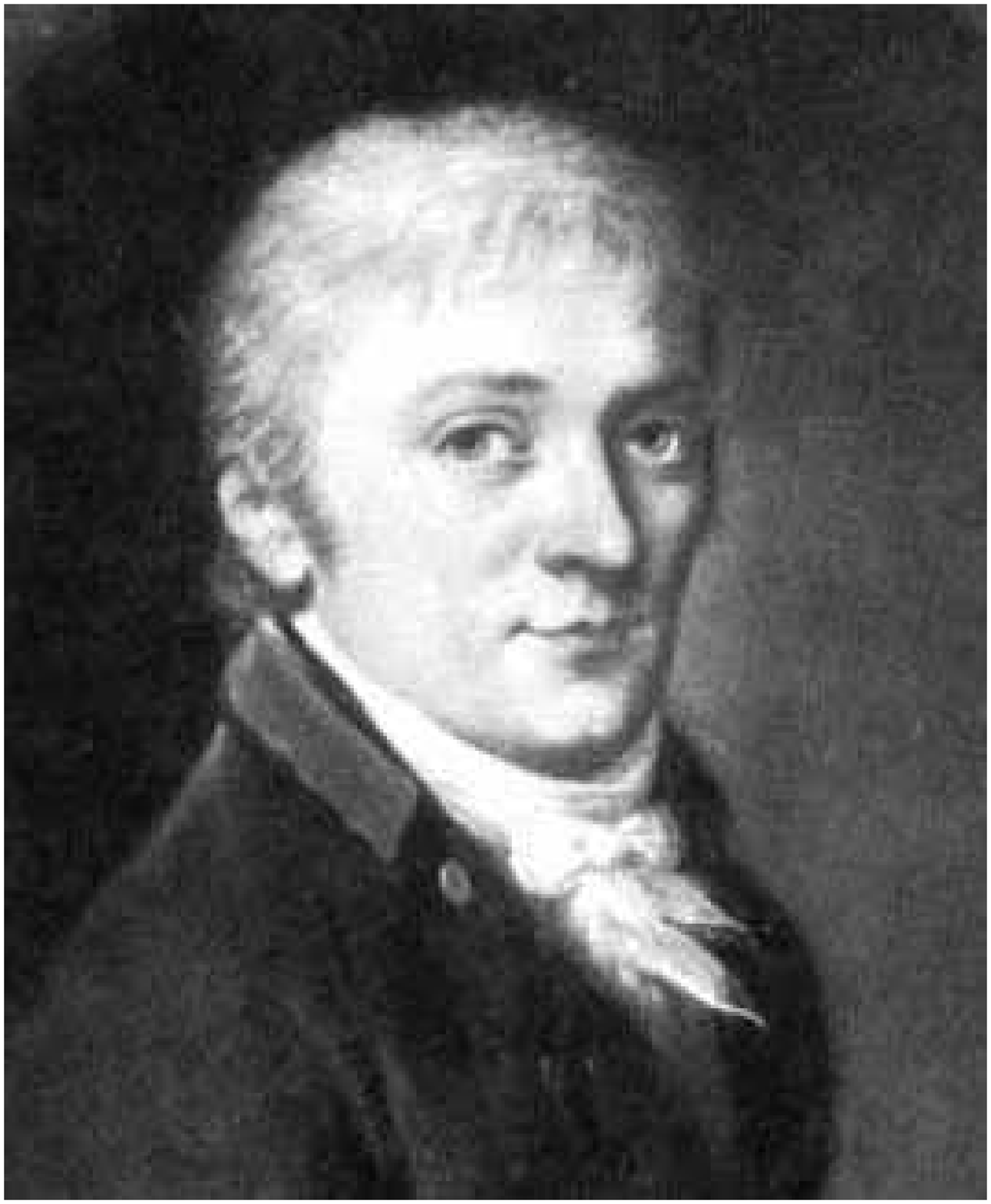}
\quad\quad
\includegraphics[width=1.5in]{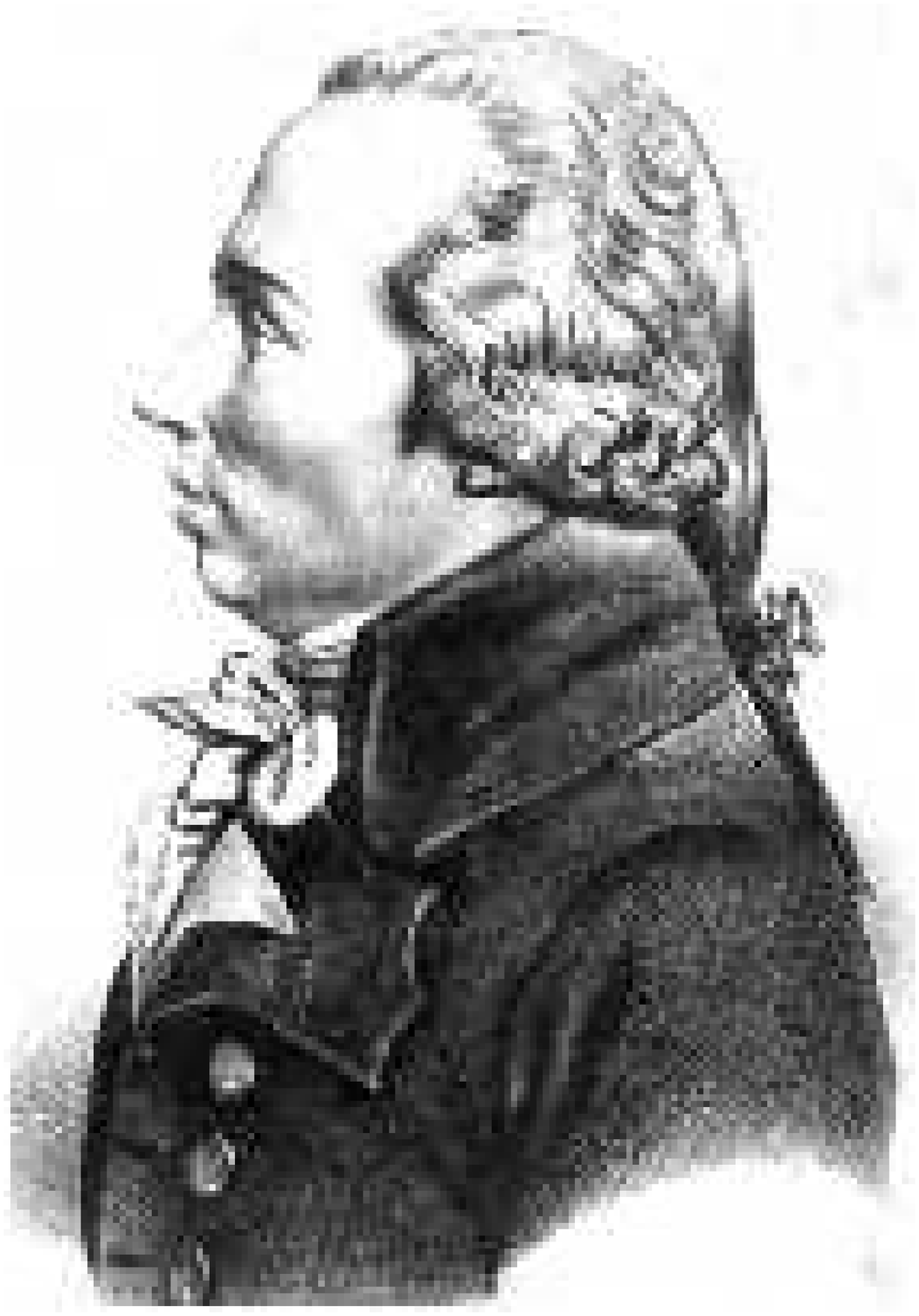}
\quad\quad
\includegraphics[width=1.5in]{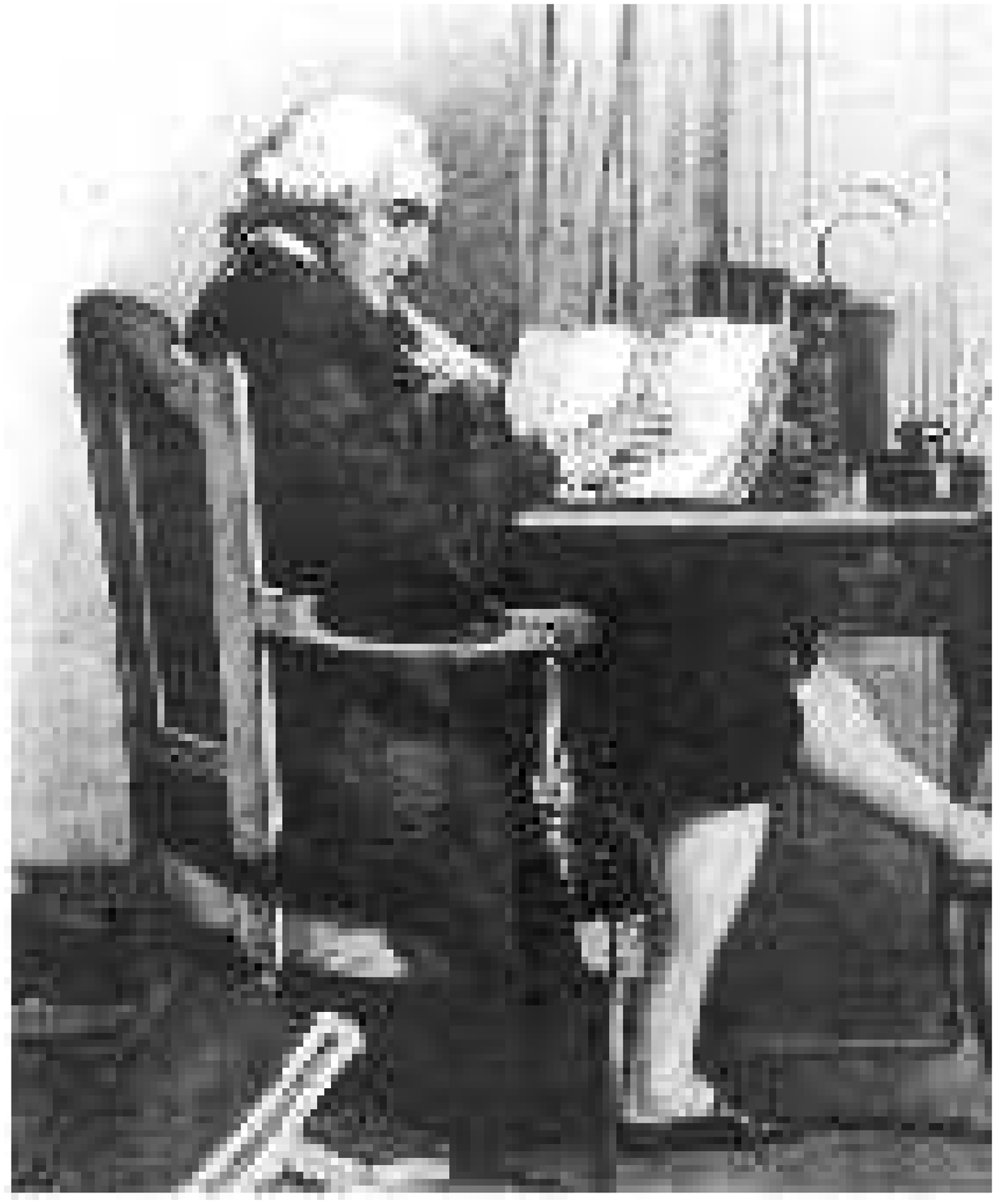}\\
F. Gauss in 1803 \qquad\qquad\quad\quad Legendre \qquad\qquad\qquad\qquad Laplace\end{center}

Secondly, Gauss himself ({\it Theoria motus coelestium}, 1809) elaborated the famous argument proving (with some 
implicit hypotheses) that once it has been assumed  the arithmetic
average is the best value to retain from among several results of quantity measurements, then,  the probability law of the error
 is necessarily the normal law. This argument has
been made more rigorous by Poincar\'e at the end of the century. Thirdly, Laplace ({Th\'eorie analytique
 des probabilit\'es}, 1811) demonstrated how the least squares method is usefull for solving linear systems when the number
of equations does not agree with the number of unknowns.
\begin{center}
\includegraphics[width=2in]{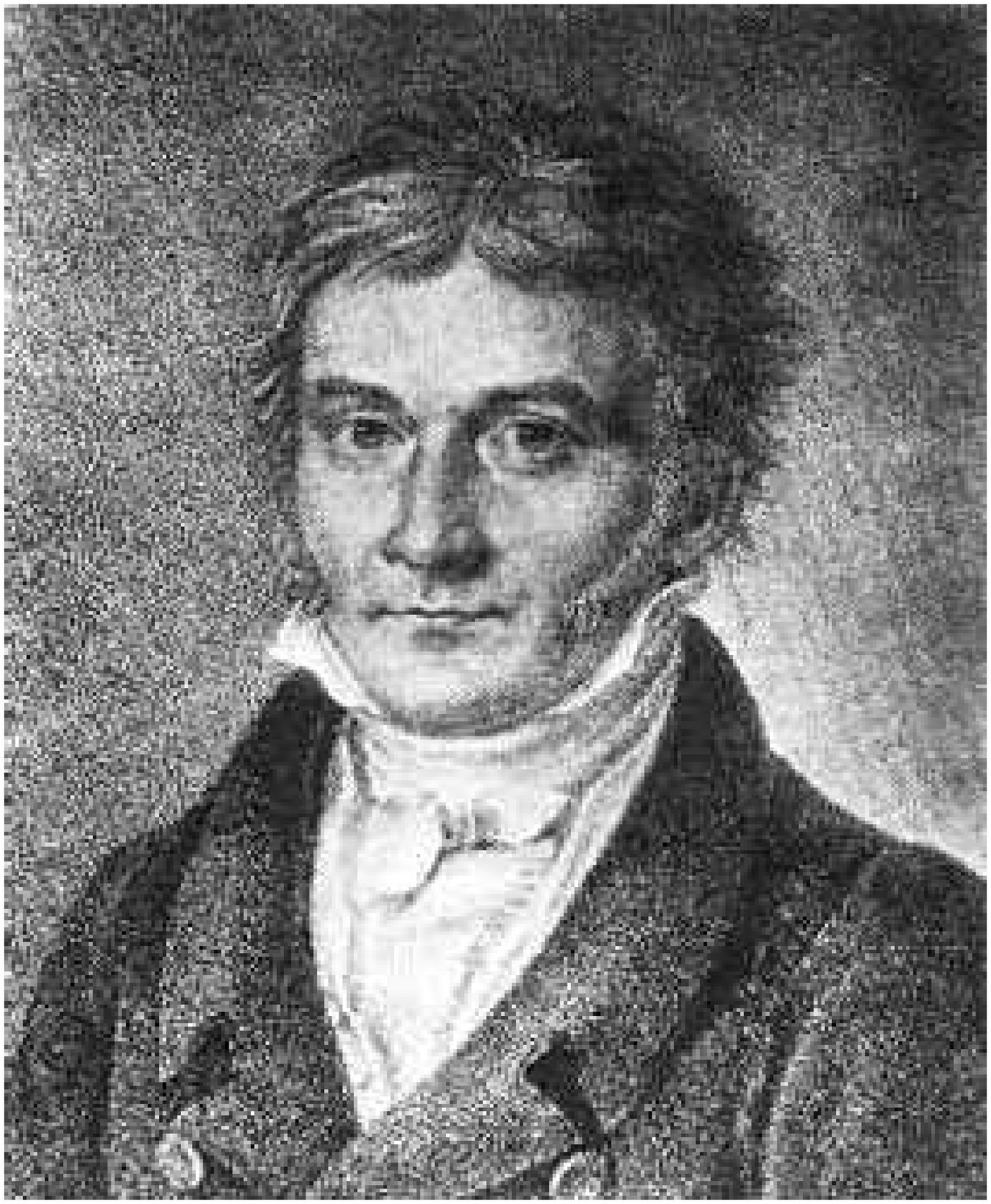}
\quad\quad
\includegraphics[width=2in]{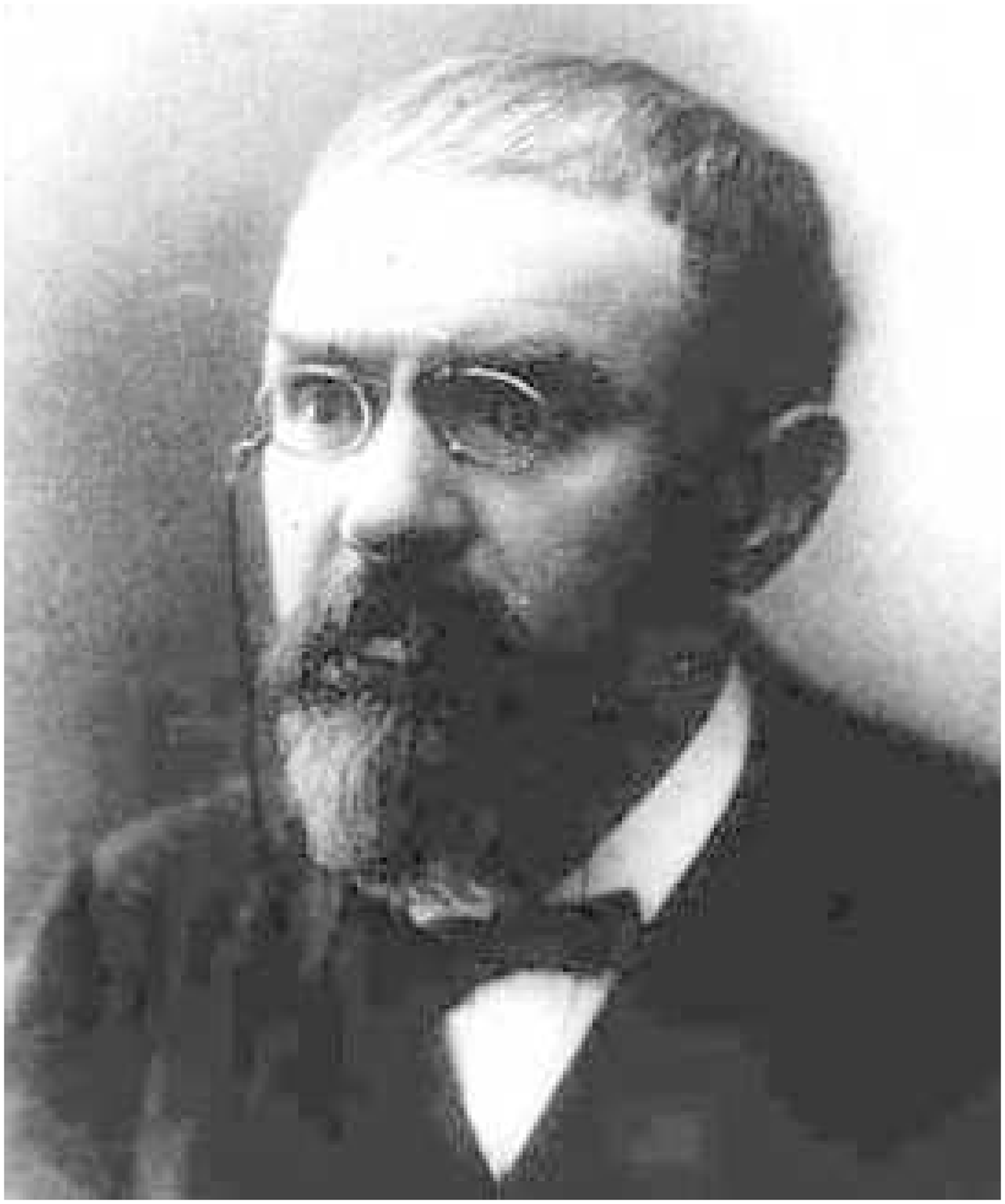}\\
F. Gauss in 1828 \qquad\qquad\qquad\quad H. Poincar\'e \end{center}
Within this same context, a few years later, Gauss became interested in the propagation of errors through calculations
 ({\it Theoria
combinationis}, 1821) and stated the following problem : 

{\it Given a quantity $U=F(V_1,V_2,V_3,\ldots)$ function of the erroneous quantities $V_1,V_2,V_3,\ldots$, compute
the potential quadratic error to expect on $U$ with the quadratic errors $\sigma_1^2,\sigma_2^2,\sigma_3^2,\ldots$
on $V_1,V_2,V_3,\ldots$ being known and assumed small and independent.}

His answer is the following formula :
\begin{equation}
\sigma_U^2=(\frac{\partial F}{\partial V_1})^2\sigma_1^2+(\frac{\partial F}{\partial V_2})^2\sigma_2^2
+(\frac{\partial F}{\partial V_3})^2\sigma_3^2+\cdots
\end{equation}
He also provides the covariance between an error on $U$ and an error on another function of the $V_i$'s.

Formula (1) displays a property which makes it much to be preferred in several respects to other
 formulae encountered in textbooks throughout the XIXth and XXth centuries. It features a coherence property. With a 
formula such as 
\begin{equation}
\sigma_U=|\frac{\partial F}{\partial V_1}|\sigma_1+|\frac{\partial F}{\partial V_2}|\sigma_1
+|\frac{\partial F}{\partial V_3}|\sigma_3+\cdots
\end{equation}
errors may depend on the way in which the function $F$ is written. Already in dimension 2, we can note that if the indentity map
were written as the composition of an injective linear map with its inverse, errors would be increased, which is hardly acceptable.

This difficulty does not arise in Gauss' calculus. Introducing the differential operator 
$$
L=\frac{1}{2}\sigma_1^2\frac{\partial^2}{\partial V_1^2}+
\frac{1}{2}\sigma_2^2\frac{\partial^2}{\partial V_2^2}+\cdots
$$
and supposing the functions to be smooth, we remark that formula (1) can be written as
$$
\sigma_U^2=L(F^2)-2FLF
$$and  coherence follows from the transport of a differential operator by an application. If $u$ and $v$ are regular
injective mappings, then, in denoting the operator $\varphi\rightarrow L(\varphi\circ u)\circ u^{-1}$ by 
$\theta_uL$, we obtain $\theta_{v\circ u}L=\theta_v(\theta_u L)$.

The errors on $V_1,V_2,V_3,\ldots$ are not necessarily supposed to be independent nor constant and may depend 
on $V_1,V_2,V_3,\ldots$ Considering a field of positive symmetric matrices $\sigma_{ij}(v_1,v_2,\ldots))$ on 
$\RR^n$ representing the conditional variances and covariances of errors on $V_1,V_2,V_3,\ldots$ given 
the values $v_1,v_2,v_3,\ldots$ of $V_1,V_2,V_3,\ldots$, then the error on $U=F(V_1,V_2,V_3,\ldots)$ 
given the values $v_1,v_2,v_3,\ldots$ of $V_1,V_2,V_3,\ldots$ is 
$$
\sigma_U^2=\sum_{ij}\frac{\partial F}{\partial V_1}(v_1,v_2,v_3,\ldots)
\frac{\partial F}{\partial V_2}(v_1,v_2,v_3,\ldots)\sigma_{ij}(v_1,v_2,v_3,\ldots)
$$
which depends solely on $F$ as mapping. This is the general form of the error calculus \`a la Gauss.\\

\noindent III.2 {\it Error propagation  through calculations : the error calculus based on Dirichlet forms}\\

The error calculus of Gauss contains the limitation of supposing that both the function $F$ and the random variables 
$V_1,V_2,V_3,\ldots$ are explicitely known. In probabilistic modelling however, we are often confronted by a situation in which
 all the random variables, functions and covariances matrices are given by limits. For such situations, a means 
of extension thereby becomes essential.

Let the quantities be defined on the probability space $(\Omega, {\cal A}, \PP)$. The quadratic error 
on a random variable $X$
is itself random, let us denote it $\Gamma[X]$. Intuitively speaking we still assume
 that the errors are infinitely small, even though this assumption does not appear in the notation. It is as though an infinitely
 small unit were
available for measuring errors fixed throughout the entire problem. The extension tool lies in the following : we assume
that if $X_n\rightarrow X$ in $L^2(\Omega, {\cal A}, \PP)$ and if the error $\Gamma[X_m-X_n]$ on $X_m-X_n$
can be made as small as we wish in $L^1(\Omega, {\cal A}, \PP)$ for $m, n$ large enough, then the error $\Gamma[X_n-X]$
on $X_n-X$ goes to zero in $L^1$.

This idea can be interpreted as a reinforced coherence principle, it means that the error on $X$ is attached to $X$ and 
furthermore, if the sequence of pairs $(X_n,{\mbox{ error on }}X_n)$ converges suitably, it converges 
necessarily to a pair $(X,{\mbox{ error on }}X)$.

The axiomatization of these idea involves the notion of closed quadratic differential
form or Dirichlet form :

{\it An error structure is a term 
$$
(\Omega, {\cal A}, \PP,\DD,\Gamma)
$$
 where $(\Omega, {\cal A}, \PP)$ is a probability space, satisfying the following properties

1) $\DD$ is a dense subvector space of $L^2(\Omega, {\cal A}, \PP)$

2) $\Gamma$ is a positive symmetric bilinear map from $\DD\times\DD$ into $L^1(\PP)$ fulfilling the functional calculus
of class ${\cal C}^1\cap{\mbox{Lip}}$, which means that if $u\in \DD^m$, $v\in\DD^n$, for $F$ and $G$ of class ${\cal C}^1$
and Lipschitz from $\RR^m$ [resp. $\RR^n$] into $\RR$, one has $F\circ u\in\DD$, $G\circ v\in\DD$ and
$$
\Gamma[F\circ u,G\circ v]=\sum_{ij}F_i^{\prime}\circ u\;\; G_i^{\prime}\circ v \;\;\Gamma[u_i,v_j]\qquad \PP-a.s.
$$
\indent 3) the bilinear form ${\cal E}[f,g]=\EE[\Gamma[f,g]]$ is closed, i.e. $\DD$ is complete under the norm
$$\| .\|_{\DD}=(\|.\|_{L^2}^2+{\cal E}[.])^{1/2}.$$}
(then the  form ${\cal E}$ is a Dirichlet form.)\\

The main benefit of the extension tool is that error theory based on Dirichlet forms extends to the
 infinite dimension, which  allows for error calculus on stochastic processes (especially on Brownian motion but also on the Poisson
space), provides several new results on stochastic differential equations, and gives applications to fluctuations in 
physics and to sensitivity analysis in finance\footnote{See the books of Malliavin, Fukushima, Ikeda-Watanabe, Bismut,
Bichteler-Gravereau-Jacod, Watanabe, Strook, Bouleau-Hirsch, Ma-R\"ockner, Nualart, \O ksendal \& al., Ustunel-Zakai,
etc. and the papers of several hundred of researchers.\\
Regarding the interpretation in terms of error propagation, see N. Bouleau, {\it Error Calculus for Finance and Physics, the
Language of Dirichlet Forms}, De Gruyter, 235p, 2003.}.\\

\noindent IV. {\large Languages with extension tools and Richard's paradox}\\

In comparing Kolmogorov's axiomatic theory of probability with the random sequences theory, we
 have emphasized for the former

- the presence of a language (syntax and semantics)

- a powerful extension tool yielding, in some sense, risky results.

This may be placed in analogy with the language of Analysis that handles real numbers.
We know, indeed, the existence of $2^{\aleph_0}$ real numbers, although only $\aleph_0$ will  ever be indicated with precision.
This is the situation highlighted by Richard's paradox (1905).
\begin{center}
\includegraphics[width=1.9in]{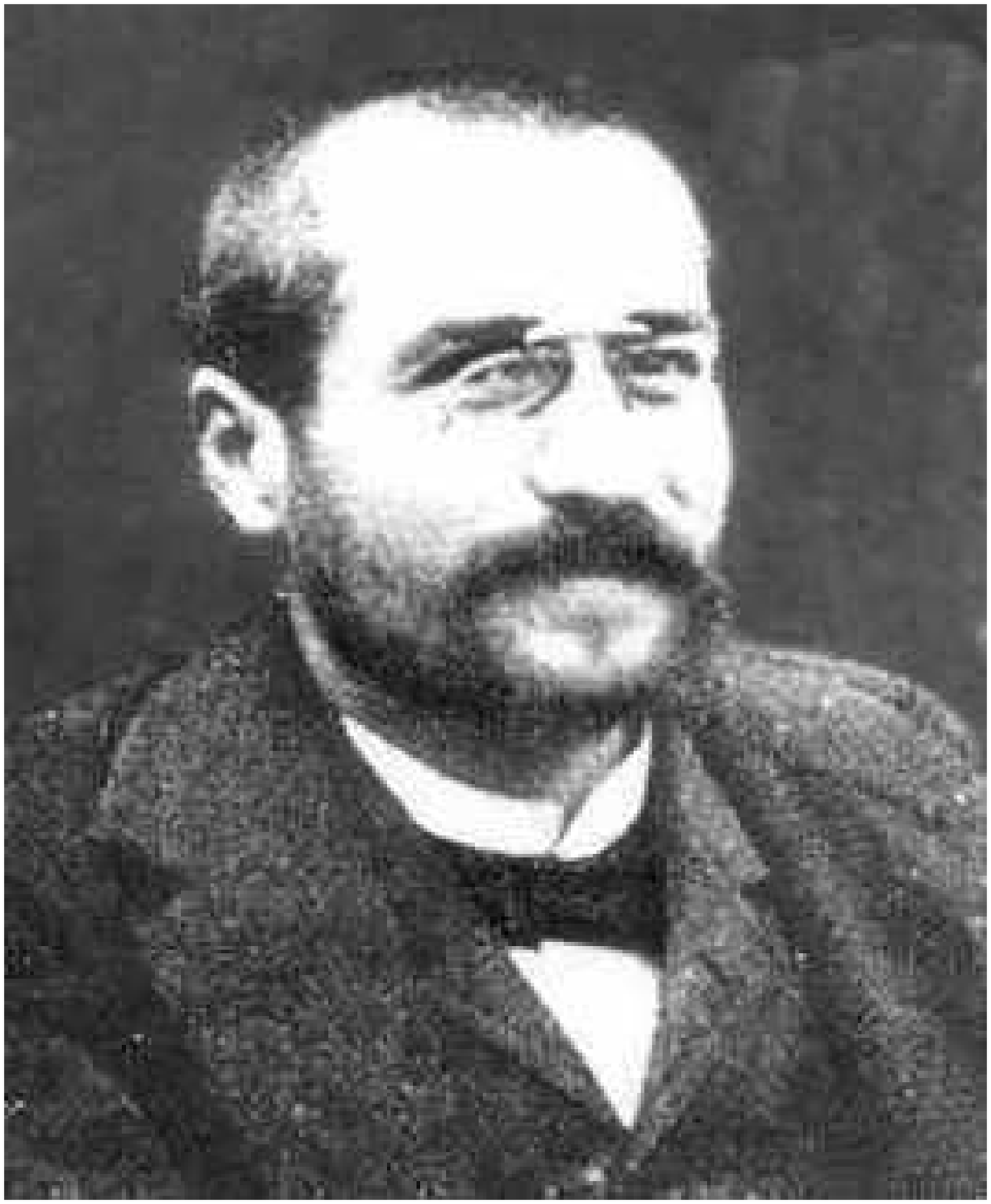}
\\
Jules Antoine Richard (1862-1956) \end{center}

The paradox can be stated as follows :

{\it Let's write all of the pairs using the 28 characters ( the 26 letters, the space and the comma to
 separate words) in alphabetic order; then the triples, and so forth, all finite sequences.
Every definition of a real number will appear in the list.

Let's cross out all the sequences which are not definitions of real numbers.

Let $u_1$ be the real number defined by the first remaining definition;

$u_2$ the one defined by the following definition;

$u_3$ the one defined by the third one;

and so forth.

We thus obtain all the real numbers defined by finitely many words, written in a particular order. The 
number {\rm a} given by the definition 
``the number without entire part, each decimal of which  immediately follows the decimal of same rank 
of the number of
same rank in the sequence $(u_n)$, the zero being considered as following the  numeral nine"
should be in the list, but cannot be equal to any number $u_n$.}

Mathematical logic is capable, of course,  of overcoming the apparent contradiction in this paradox. Nevertheless, a true
phenomenon has indeed been highlighted : there are $2^{\aleph_0}$ real numbers, we dont know how large this cardinal 
$2^{\aleph_0}$ actually is, and only $\aleph_0$ real numbers will  ever be precisely defined.

In such a situation,  for mathematical Analysis, we have chosen a language with an extension tool : the Cauchy criterion. This strategy
allows handling real numbers defined by limits regardless of the  construction of the used convergent sequence. This tool has then been carried
from the real case to the functional case by the notions of Hilbert space and Banach space which are certainly
ones of the most powerful concepts of XXth century Analysis.

\end{document}